\documentclass[a4paper]{amsart}

\usepackage{graphicx}
\usepackage{float}
\usepackage{url}
\usepackage{enumerate}
\usepackage{xspace}
\usepackage{multirow}
\usepackage{hhline}
\usepackage{microtype}
\usepackage{lmodern}
\usepackage[font=small, labelfont=bf]{caption}
\usepackage{soul}
\usepackage{subcaption}

\newcommand{\df}[1]{\ensuremath{\mathrm{df}({#1})}\xspace}
\newcommand{\rdf}[1]{\ensuremath{\mathrm{rdf}({#1})}\xspace}

\newcommand{\core}{{\mathrm{core}}}

\theoremstyle{plain}
\newtheorem{theorem}{Theorem}[section]
\newtheorem{lemma}[theorem]{Lemma}
\newtheorem{proposition}[theorem]{Proposition}
\newtheorem{corollary}[theorem]{Corollary}

\theoremstyle{definition}
\newtheorem{definition}[theorem]{Definition}

\newtheorem*{example*}{Example}
\newtheorem{problem}{Problem}

\title[Regular colouring defect of a cubic graph]
{Regular colouring defect of a cubic graph and
the conjectures of Fan-Raspaud and Fulkerson}

\author[J. Karab\'a\v{s}]{J\'an Karab\'a\v{s}}
\email[J. Karab\'a\v{s}]{jan.karabas@stuba.sk}
\author[E. M\'a\v cajov\'a]{Edita M\'a\v cajov\'a}
\email[E. M\'a\v cajov\'a]{macajova@dcs.fmph.uniba.sk}
\author[R. Nedela]{Roman Nedela}
\email[R. Nedela]{nedela@savbb.sk}
\author[M. \v Skoviera]{Martin \v Skoviera}
\email[M. \v Skoviera]{skoviera@dcs.fmph.uniba.sk}

\address[J. Karab\'a\v{s}]{
Faculty of Electrical Engineering and Information Technology, Slovak University of Technology, Bratislava, Slovakia
}

\address[E. M\'a\v cajov\'a, M. \v Skoviera]{
Comenius University, Mlynsk\' a dolina, Bratislava, Slovakia
}

\address[R. Nedela]{
Faculty of Applied Sciences, University of West Bohemia,
Pilsen, Czech Republic}

\address[J. Karab\'a\v{s}, R. Nedela]{
Mathematical Institute of Slovak Academy of Sciences,
Bansk\'a Bystrica, Slovakia
}

\begin{document}

\begin{abstract}
We introduce a new invariant of a cubic graph -- its regular
colouring defect -- which is defined as the smallest number of edges left uncovered by any collection of three perfect matchings that
have no edge in common. This invariant is a modification of
colouring defect, an invariant introduced by Steffen (J. Graph Theory 78
(2015), 195--206), whose definition does not require the empty
intersection condition. In this paper we discuss the
relationship of this invariant to the well-known conjectures of
Fulkerson (1971) and Fan and Raspaud (1994) and prove that
colouring defect and regular colouring defect can be arbitrarily far apart.
\end{abstract}
\maketitle

\section{Introduction}\label{sec:intro}
\noindent{}In this paper we introduce a new structural invariant
of cubic graphs, called regular colouring defect. The \textit{regular colouring} defect of a bridgeless cubic graph $G$, denoted by $\rdf{G}$,  is defined as the smallest number of edges of $G$
left uncovered by three perfect matchings $M_1$, $M_2$, and
$M_3$ such that $M_1\cap M_2\cap M_3=\emptyset$. Since
$\rdf{G}=0$ if and only if $G$ is $3$-edge-colourable, regular
colouring defect can be taken as one of many measures of uncolourability,
widely discussed in a recent survey article by Fiol et al.
\cite{FMS-survey}.

Regular colouring defect is closely related to an invariant introduced
by Steffen in 2015 in \cite{S2} as $\mu_3$ and later termed
\textit{colouring} defect which is simply the smallest number of edges left uncovered by
any three perfect matchings. The latter invariant was
extensively studied in \cite{KMNS-girth, KMNS-Berge, KMNS-red,
KMNS-eurocomb, KMNS-Bergegen}. It is clear from the definition
that colouring defect never exceeds regular colouring defect.

The primary motivation for the study of regular colouring defect is its
close relationship to the conjectures of Fulkerson (1971) and
Fan and Raspaud (1994); see \cite{F} and \cite{FR},
respectively. The former conjecture suggests that every
bridgeless cubic graph admits a collection of six perfect
matchings such that each edge belongs to exactly two of them,
while the latter conjecture states that there exist three
perfect matchings with empty intersection. The conjecture of
Fan and Raspaud, if true, implies that three perfect matchings required by the definition of a regular colouring defect indeed exist in every bridgeless cubic graph. Another
reason why regular colouring defect deserves particular attention is the
fact that every triple of perfect matchings $M_1$, $M_2$, and
$M_3$ with $M_1\cap M_2\cap M_3=\emptyset$ induces a
nowhere-zero flow whose values can be taken as points of the
Fano plane and whose flow patterns around vertices correspond
to lines of the Fano plane. Hence, regular colouring defect can be
investigated with the help of techniques developed in the
theory of nowhere-zero flows. These ideas will be discussed in
detail in Section~\ref{sec:fulkerson}.

The main purpose of this paper is to show that the difference
between colouring defect and regular colouring defect can be arbitrarily large.

\begin{theorem}\label{thm:main}
For every integer $d\ge 4$ there exists a bridgeless cubic
graph $G$ whose colouring defect equals $4$ and regular colouring defect is at least $d$.
\end{theorem}

This result, which will be proved in Section~\ref{sec:main}, does not extend
to $d=3$, because in the latter case colouring defect and regular colouring defect of a cubic graph coincide (see Corollary~\ref{cor:rdf-df}).

\section{Preliminaries}\label{sec:prelim}
\noindent{}Graphs considered in this paper are finite and
undirected; parallel edges and loops are permitted. A
\emph{circuit} is a connected $2$-regular graph. A
$k$-\emph{cycle} is  a circuit of length~$k$. The length of a
shortest circuit in a graph is its \emph{girth}. A graph $G$ is
said to be \emph{cyclically $k$-edge-connected} if the removal
of fewer than $k$ edges from $G$ cannot create a graph with at
least two components containing circuits. An edge cut $S$ in
$G$ that separates two circuits from each other is
\emph{cycle-separating}.

Large cubic graphs can be conveniently described  as
combinations of smaller building blocks called multipoles,
which, in contrast to graphs, may contain dangling edges or
even isolated edges. Each \emph{multipole} $M$ has its vertex
set $V(M)$, its edge set $E(M)$, and an incidence relation
between vertices and edges. Each edge of $M$ has two
\emph{ends}, and each end may, but need not be, incident with a
vertex of $M$. An end of an edge that is not incident with a
vertex is called a \emph{free end} or a \emph{semiedge}. An
edge with exactly one free end is called a \emph{dangling
edge}. An \emph{isolated edge} is an edge whose both ends are
free. All multipoles considered in this paper are \emph{cubic};
it means that every vertex is incident with exactly three
edge-ends. An \emph{$n$-pole} $M$ is a multipole with $n$ free
ends. Free ends of a multipole can be distributed into pairwise
disjoint sets, called \emph{connectors}.  An \emph{$(n_1,
n_2,\ldots, n_k)$-pole} is an $n$-pole with $n=n_1+n_2+\cdots +
n_k$ whose semiedges are distributed into $k$ connectors
$S_1,S_2,\ldots, S_k$, each $S_i$ being of size $n_i$. A
\emph{dipole} is a multipole with two connectors, while a
\emph{tripole} is a multipole with three connectors.

A \emph{$3$-edge-colouring} of a cubic graph or a cubic
multipole is an assignment of three colours (usually $1$, $2$,
and $3$) to the edges of $G$ such that any two adjacent edges
receive distinct colours. A $2$-connected cubic graph that does
not admit a $3$-edge-colouring is called a \emph{snark}.

It is
often convenient to regard the colours $1$, $2$, and $3$ as
nonzero elements of the group
$\mathbb{Z}_2\times\mathbb{Z}_2$, say $1=(0,1)$,
$2=(1,0)$, and $3=(1,1)$. The condition requiring the three
colours meeting at any vertex $v$  be distinct becomes
equivalent to requiring that the sum of colours at $v$ be zero.
In other words, a $3$-edge-colouring of a cubic graph is a
nowhere-zero $\mathbb{Z}_2\times\mathbb{Z}_2$-flow and the
condition on the sum of colours is nothing but Kirchhoff's law.
To be more precise, given an abelian group $A$, a
\emph{nowhere-zero $A$-flow} on a multipole $M$ is an
assignment of an orientation to the edges of $M$ and a function
$\phi\colon E(M)\to A-\{0\}$ such that at each vertex the sum
of incoming values equals the sum of outgoing values. Note that
the orientation can be ignored whenever every element of
$A-\{0\}$ has order~$2$, in which case $A$ is isomorphic to $\mathbb{Z}_2^n$ for
some $n\ge 1$. If $A=\mathbb{Z}$ and there is an
integer $k\ge 2$ such that $|\phi(e)|<k$, we speak of a
\emph{nowhere-zero $k$-flow}. It is well known that a graph
admits a nowhere-zero $k$-flow if and only if it admits a
nowhere-zero $A$-flow where $A$ is an arbitrary abelian group
of order~$k$ (see Diestel \cite[Corollary~6.3.2 and
Theorem~6.3.3]{D}). This result can easily be extended to
multipoles.

\medskip

The following well-known statement is a direct consequence of
Kirchhoff's law. It tells us that the total outflow from any
nonempty set of vertices equals $0$.

\begin{lemma}[Parity Lemma]\label{lem:par}
Let $M$ be a multipole with $n$ free edge-ends $s_1$, $s_2$,
\ldots, $s_n$. If $\sigma$ is an arbitrary $3$-edge-colouring
of $M$, then
$$\sum_{i=1}^n\sigma(s_i)=0.$$
Equivalently, the number of free ends of $M$ carrying any fixed colour has the same parity as $n$.
\end{lemma}

\section{Regular colouring defect of a cubic graph}\label{sec:fulkerson}
\noindent{}Every $3$-edge-colouring of a cubic graph $G$
determines three pairwise disjoint perfect matchings, the
colour classes. Any collection $\mathcal{M}=\{M_1,
M_2,\penalty0 M_3\}$ of three perfect matchings of $G$ can
therefore be regarded as a generalisation of a
$3$-edge-colouring. We call such a collection  a
\emph{$3$-array} of perfect matchings of $G$. An edge of $G$
that belongs to at least one of the perfect matchings from
$\mathcal{M}$ is said to be \emph{covered}. An edge is
\emph{uncovered}, \emph{simply covered}, \emph{doubly covered},
or \emph{triply covered} if it belongs, respectively, to zero,
one, two, or three distinct members of~$\mathcal{M}$.

If a cubic graph $G$ is $3$-edge-coloured, then all its edges
are simply covered. If $G$ cannot be  $3$-edge-coloured, then
at least one edge remains uncovered under any $3$-array of
perfect matchings. The minimum number of uncovered edges taken
over all $3$-arrays for $G$ is the \emph{colouring defect}
(briefly, the \emph{defect}) of~$G$; it is denoted by
$\df{G}$. Since $\df{G}=0$ if and only if $G$ is
$3$-edge-colourable, colouring defect can measure to what extent
cubic graphs differ from $3$-edge-colourable ones.

Given an arbitrary $3$-array $\mathcal{M}=\{M_1, M_2, M_3\}$ of
a cubic graph $G$ we define the \emph{core} of $\mathcal{M}$,
denoted by $\core(\mathcal{M})$, to be the subgraph induced by
the set of edges of $G$ that are not simply covered. A
\emph{$k$-core} is a core with $k$ uncovered edges. A $k$-core,
as well as the corresponding $3$-array, will be called
\emph{optimal} if $k=\df{G}$. It may be worth mentioning that
if $G$ is $3$-edge-colourable and $\mathcal{M}$ consists of
three pairwise disjoint perfect matchings, then $\mathcal{M}$
is optimal and $\core(\mathcal{M})$ is empty. On the other hand, if $G$ has no perfect matching, then each edge of $G$ is uncovered for every $3$-array $\mathcal{M}$ of $G$, and hence $\core(\mathcal{M})=G$. In particular, the colouring defect of $G$ equals the number of edges.

The following proposition, which extends Lemma~2.2 and
Corollary~2.5 of \cite{S2}, lists the most fundamental
properties of cores.

\begin{proposition}[{\cite[Propositon~3.2]{KMNS-red}}]\label{prop:core}
Let $\mathcal{M}=\{M_1, M_2, M_3\}$ be an arbitrary $3$-array
of perfect matchings of a snark $G$. Then the following hold:
\begin{enumerate}[{\rm (i)}]
\item Every component of $\core(\mathcal{M})$ is
    either an even circuit, which alternates doubly covered
    and uncovered edges, or a subdivision of a cubic graph.
    Moreover, $\core(\mathcal{M})$ is a collection of
    disjoint circuits if and only if $G$ has no triply
    covered edge.
\item Every $2$-valent vertex of $\core(\mathcal{M})$ is
    incident with one doubly covered edge and with one uncovered
    edge, while every $3$-valent vertex is incident with
    one triply covered edge and two uncovered edges.
\item Every snark $G$ has $\df{G}\ge 3$. If $\df{G}=3$,
    then every optimal core is an induced $6$-cycle.
\end{enumerate}
\end{proposition}

We say that a $3$-array $\mathcal{M}$ for a cubic graph $G$ is
\emph{regular} if $\core(\mathcal{M})$ is empty, or each component of $\core(\mathcal{M})$ is a
circuit, that is, if the core is a $2$-regular subgraph of $G$.
By Proposition~\ref{prop:core}(ii), a $3$-array $\mathcal{M}$
is regular if and only if it has no triply covered edges. We
can now define the \emph{regular colouring} defect of~$G$, denoted by
$\rdf{G}$, to be the minimum number of uncovered edges taken
over all regular $3$-arrays for $G$. 
In particular, if $G$ has no regular $3$-array, then $\rdf{G}$ equals the number of edges of $G$. A regular $3$-array $\mathcal{M}$ for which the number of uncovered edges equals
$\rdf{G}$ will be called \emph{regularly optimal}. It is clear
from the definition that for every bridgeless cubic graph $G$
one has
$$\rdf{G}\ge\df{G}.$$

As previously mentioned, Fan and Raspaud \cite{FR} made a
conjecture that every bridgeless cubic graph has three perfect
matchings $M_1$, $M_2$, and $M_3$ with $M_1\cap M_2\cap
M_3=\emptyset$. In other words, they conjectured that every
bridgeless cubic graph admits a regular $3$-array.
This conjecture remains widely open, nonetheless,
it is known \cite{MS-Comb} that the conjecture holds for all
bridgeless cubic graphs with oddness at most $2$ (that is,
graphs that contain a $2$-factor with at most two odd
circuits).

The following statement is a corollary of Proposition~\ref{prop:core}.

\begin{corollary}\label{cor:rdf-df}
For an arbitrary snark $G$, $\rdf{G}=3$ if and only if $\df{G}=3$.
\end{corollary}

\begin{proof}
If $G$ is a snark with $\rdf{G}=3$, then, by
Proposition~\ref{prop:core}~(iii) we have
$3\le\df{G}\le\rdf{G}=3$, so $\df{G}=3$. Conversely, assume
that $\df{G}=3$, and let $\mathcal{M}=\{M_1,M_2,M_3\}$ be an
optimal $3$-array for $G$. By Proposition~\ref{prop:core}
again, the core of $\mathcal{M}$ is a $6$-cycle that alternates
uncovered and doubly covered edges. It follows that
$\mathcal{M}$ is a regular $3$-array, whence $\rdf{G}=3$.
\end{proof}

\begin{figure}[h!]
	\includegraphics[scale=1.4]{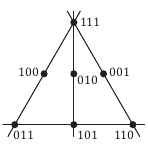}\label{fig:konf2}
	\caption{Four lines of the Fano plane that correspond to regular $3$-arrays}
	\label{fig:konf}
\end{figure}

Given a regular $3$-array $\mathcal{M}=\{M_1, M_2, M_3\}$ for
$G$, we can define a mapping
\[
\chi\colon E(G)\to \mathbb{Z}_2^3,
\quad e\mapsto \chi(e)=(x_1, x_2, x_3)
\]
by setting $x_i = 0$ if and only if $e\in M_i$ for
$i\in\{1,2,3\}$. Since the complement of each $M_i$ in $G$ is a
$2$-factor, it is easy to see that $\chi$ is a
$\mathbb{Z}_2^3$-flow. We call $\chi$ \emph{the characteristic
flow} for $\mathcal{M}$. The fact that $\mathcal{M}$ is regular
implies that the value $(0,0,0)$ never occurs, which means that
$\chi$ is, in fact, a nowhere-zero flow. Consequently, we can
regard the flow values of $\chi$ as points of the Fano plane
$PG(2,2)$ represented by the standard projective coordinates
from $\mathbb{Z}_2^3-\{0\}$. The definition of the
characteristic flow further implies that for every vertex $v$
of $G$ the three flow values meeting at $v$ form a line of the
Fano plane, that is, a triple $\{x,y,z\}$ of points such that
$x+y+z=0$. Since every vertex of $G$ is incident with exactly
one edge of each member of $\mathcal{M}$, one can easily deduce
that only four lines of the Fano plane can occur around a
vertex, and these lines form a configuration shown in
Figure~\ref{fig:konf}.
Therefore, the characteristic flow of a
regular array is a nowhere-zero flow with an additional
geometric structure. It is also worth mentioning that the
characteristic flow of $\mathcal{M}$ is at the same time a
proper edge colouring. Further investigation of regular defect
thus can employ both methods of nowhere-zero flows as well as
those of edge colourings. These methods have been already
applied in \cite{KMNS-red} for the study of snarks whose defect
-- and by Corollary~\ref{cor:rdf-df} also regular defect --
equals~3.

So far we have discussed regular $3$-arrays and their relation to the Fan-Raspaud conjecture. We now proceed to examining their relation to the Fulkerson conjecture, which states that
\emph{every bridgeless cubic graph $G$ admits a collection of
six perfect matchings (not necessarily distinct) such that each
edge belongs to precisely two of them}. Such a collection will
be called a \emph{Fulkerson cover} of $G$. Our discussion features the following important concept.

\begin{definition}
Two regular $3$-arrays $\mathcal{M}_1$ and
$\mathcal{M}_2$ for a cubic graph $G$ are said to be
\emph{complementary} if
$\core(\mathcal{M}_1)=\core(\mathcal{M}_2)$ but the
corresponding sets of uncovered edges are disjoint.
\end{definition}

The following result explains the significance of pairs of complementary $3$-arrays; it includes Lemma~3.1 of
\cite{Hao+CQ_BF} and Theorem~1.1 of \cite{Chen+Fan}.

\begin{theorem}\label{thm:complementary}
The following three statements are equivalent for every
bridgeless cubic graph $G$.
\begin{enumerate}[{\rm (i)}]
\item $G$ has a Fulkerson cover.
\item $G$ contains a set of disjoint even circuits which
    constitute the common core of two complementary regular
    $3$-arrays of perfect matchings.
\item $G$ contains two disjoint matchings $P_1$ and $P_2$
    such that $P_1\cup P_2$ is a set of circuits of $G$ and
    each of $G-P_1$ and $G-P_2$ admits a nowhere-zero
    $4$-flow.
\end{enumerate}
\end{theorem}

\begin{proof}
(i) $\Rightarrow$ (ii): Let
$\mathcal{F}=\{M_1,M_2,M_3,M_4,M_5,M_6\}$ be a Fulkerson cover
of $G$. Consider a $3$-array $\mathcal{M}=\{M_1,M_2,M_3\}$.
Since $\mathcal{M}\subseteq\mathcal{F}$, the core of
$\mathcal{M}$ does not contain any triply covered edge. It
follows that $\core(\mathcal{M})$ is regular, that is, a
disjoint union of even circuits. The sets of simply
covered edges for $3$-arrays $\mathcal M=\{M_1,M_2,M_3\}$ and
$\mathcal{M}'=\{M_4,M_5,M_6\}$ coincide, therefore the $3$-arrays
$\mathcal{M}$ and $\mathcal{M}'$ have the same core with the
uncovered and doubly covered edges interchanged. Thus,
$\mathcal{M}$ and $\mathcal{M}'$ are complementary regular
$3$-arrays for $G$, as required.

(ii) $\Rightarrow$ (iii):  Assume that we have a pair of
complementary regular $3$-arrays $\mathcal{M}_1$ and
$\mathcal{M}_2$ for~$G$ with common core $Q$. For each
$i\in\{1,2\}$ let $P_i$ denote the matching consisting the
edges left uncovered by the $3$-array $\mathcal{M}_i$. The
definition of complementary $3$-arrays yields that $P_1\cup
P_2=Q$, so $P_1\cup P_2$ is a set of circuits of $G$. We now
show that each $G-P_i$ admits a nowhere-zero $4$-flow. Due to
the well-known equivalence results mentioned earlier it is
sufficient to find a nowhere-zero
$\mathbb{Z}_2\times\mathbb{Z}_2$-flow $\phi_i$ on $G-P_i$.
Assume that $\mathcal{M}_i=\{M^i_1,M^i_2,M^i_3\}$ and let an
edge $x\in E(G-P_i)$ be arbitrary. Since each edge of $G-E(Q)$
is simply covered by $\mathcal{M}_i$, for such $x$ we can
unambiguously set $\phi_i(x)=t\in\{1,2,3\}$ whenever $x\in
M^i_t$. It remains to determine $\phi_i(x)$ when $x$ belongs to
$Q-P_i$, that is, when $x$ is doubly covered by
$\mathcal{M}_i$. Assume that $x\in M^i_j\cap  M^i_k$ with $i\ne
j$. Then both edges of $G-P_i$ adjacent to $x$ belong to the
third perfect matching $M^i_l$, where $\{j,k,l\}=\{1,2,3\}$,
and we can set $\phi_i(x)=l$. The resulting mapping
$\phi_i\colon
E(G-P_i)\to\{1,2,3\}=\mathbb{Z}_2\times\mathbb{Z}_2-\{0\}$ is
easily seen to fulfil the Kirchhoff law, so $\phi_i$ is a
nowhere-zero $\mathbb{Z}_2\times\mathbb{Z}_2$-flow on $G-P_i$.

(iii) $\Rightarrow$ (i): Assume that $G$ contains two disjoint
matchings $P_1$ and $P_2$ such that $P_1\cup P_2$ is a set
circuits of $G$ and both $G-P_1$ and $G-P_2$ admit a
nowhere-zero $4$-flow. We construct a Fulkerson cover for $G$
by specifying a suitable proper edge
colouring $\xi\colon
E(G)\to\mathcal{P}_{6,2}$, where $\mathcal{P}_{6,2}$ denotes
the set of all $2$-element subsets of the set
$\{1,2,3,4,5,6\}$. Recall that for each $i\in\{1,2\}$ we have a
nowhere-zero flow $\phi_i\colon
E(G-P_i)\to\{1,2,3\}=\mathbb{Z}_2\times\mathbb{Z}_2-\{0\}$. We
now use these two flows to define the colouring $\xi$, this
time interpreting the values $1$,~$2$, and $3$ as integers. For
each edge $x$ of $G-(P_1\cup P_2)$ we set
$\xi(x)=\{\phi_1(x),\phi_2(x)+3\}$, for $x\in P_1$ we set
$\xi(x)=\{1,2,3\}-\{\phi_1(x)\}$, and for $x\in P_2$ we set
$\xi(x)=\{4,5,6\}-\{\phi_2(x)+3\}$ (addition in~$\mathbb{Z}$).
One can easily check that at each vertex $v$ of $G$ the colours
of the three edges incident with $v$ form a partition of the
set $\{1,2,3,4,5,6\}$ into $2$-element subsets. It follows that
if we set $M_i=\{e\in E(G);\, i\in\xi(e)\}$, then for each
$i\in \{1,2,3,4,5,6\}$ the set $M_i$ is a perfect matching and
the collection $\{M_1,M_2,M_3,M_4,M_5,M_6\}$ is a Fulkerson
cover of $G$.
\end{proof}

\section{Proof of Theorem~\ref{thm:main}}\label{sec:main}
We begin with the following easy but useful proposition.

\begin{proposition}\label{prop:rdf-girth}
If $G$ is a snark of girth $g$, then $\rdf{G}\ge g/2$.
\end{proposition}

\begin{proof}
Take a regularly optimal $3$-array
$\mathcal{M}=\{M_1,M_2,M_3\}$ for $G$. By
Proposition~\ref{prop:core}~(i), its core is a collection of
disjoint even circuits that alternate uncovered and doubly
covered edges. Let $C$ be a $k$-cycle constituting a
component of $\core(G)$. Then
$$\rdf{G}\ge k/2\ge g/2,$$
as required. 
\end{proof}

Now we can proceed to the main results of this paper. A pair
$\{u,v\}$ of vertices of a snark $G$ is said to be
\emph{non-removable} if $G-\{u,v\}$ is $3$-edge-colourable.

\begin{theorem}\label{thm:2triangles}
Let $G$ be a snark with  a non-removable pair of adjacent
vertices $\{u,v\}$, and let $\bar G$ be the snark arising
from $G$ by inflating both $u$ and $v$ to a triangle. Then
$3\le \df{\bar G}\le 4$.
\end{theorem}

\begin{figure}[h!]
	\includegraphics[width=0.4\textwidth]{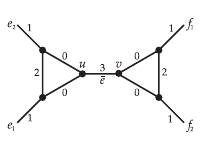}
	\caption{An irregular $4$-core with edge-labels showing the
		number of perfect matchings containing the edge in question}
	\label{fig:2triangles}
\end{figure}

\begin{proof}
Clearly, the graph $\bar G$ is again a snark, so $\df{\bar
G}\ge 3$, by Proposition~\ref{prop:core}. We prove that
$\df{\bar G}\le 4$. Let $e_1$ and $e_2$ be the edges of $G$
incident with $u$ and different from $e=uv$, and similarly let
$f_1$ and $f_2$ be the edges incident with $v$ and different
from $uv$. Take an arbitrary  $3$-edge-colouring $\sigma$ of
$G-\{u,v\}$. Parity lemma implies that if we keep the four
dangling edges, then $\sigma(e_1)=\sigma(e_2)$ and
$\sigma(f_1)=\sigma(f_2)$, otherwise $\sigma$ could be extended
to a $3$-edge-colouring of the entire graph $G$. Now consider the
graph $\bar G$. Clearly, for each edge $x$ of $G$ there is a
unique corresponding edge $\bar x$ in $\bar G$. Consequently,
the three colour classes of the colouring $\sigma$ of $G$
extend to a $3$-array of $\bar G$ where the edge $\bar e$ is
triply covered, the four edges adjacent to $\bar e$ are
uncovered, the third edge of each triangle is doubly covered,
and all remaining edges of $\bar G$ are simply covered, see
Figure~\ref{fig:2triangles}. The corresponding core is an
irregular $4$-core. Consequently, $\df{\bar G}\le 4$.
\end{proof}

The following theorem was proved in \cite{KMNS-girth}.

\begin{theorem}\label{thm:largegirth}
For every even integer $g\ge 6$ there exists a cyclically
$5$-edge-connected snark $G_g$ of girth $g$ that contains a
non-removable pair of adjacent vertices.
\end{theorem}

We sketch the construction of the family $(G_g)_{g\ge 6}$, which
constitutes the core of the proof of
Theorem~\ref{thm:largegirth}. We also indicate why $G_g$ has
all the stated properties. For details we refer the reader to \cite{KMNS-girth}.

The graphs $G_g$ are constructed by employing the method of
superposition. Its main idea is to `inflate' a given base graph
$G$ into a large cubic graph $\tilde G$ by substituting
vertices of $G$ with `fat vertices' (tripoles), called
\emph{supervertices}, and edges of $G$ with `fat edges'
(dipoles), called \emph{superedges}. For formal definitions and
a detailed description of the superposition method we refer the
reader to \cite{Ko, MS-superp}.

\begin{figure}[h!]
\centering
\hfill\begin{subfigure}[c]{0.35\textwidth}
\includegraphics[width=\textwidth]{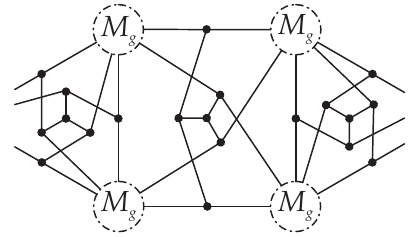}
\end{subfigure}\hfill
\begin{subfigure}[c]{0.35\textwidth}
\includegraphics[width=\textwidth]{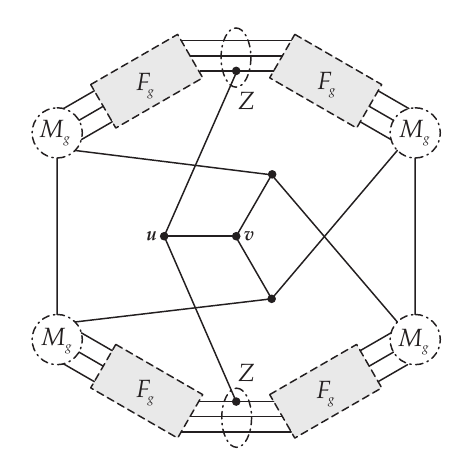}
\end{subfigure}\hfill\phantom{x}
\caption{The superedge $F_g$ (left) and the resulting snark $G_g$ (right)}
\label{fig:superpos}
\end{figure}

The base graph for our construction of $G_g$ is the Petersen
graph $Pg$. In order to create the desired snarks from $Pg$ we
need two types of supervertices, represented by a
$(3,3,1)$-pole $Z$ and a $(3,1,1)$-pole $M_g$, and one type of
a superedge, represented by a $(3,3)$-pole~$F_g$. 

We define $Z$
to be the unique trivalent $(3,3,1)$-pole of order 1; it
consists of a vertex incident with three semiedges contributing
to three different connectors and of two isolated edges each of
which contributes to different connectors of size $3$, see
Figure~\ref{fig:superpos} on the right-hand side.
The $(3,1,1)$-pole $M_g$ is created
from a connected bipartite vertex-transitive cubic graph $L_g$
of girth $g$ by removing a path of length $2$. Such a graph
indeed exists: Theorem~4.8 in \cite{NS-largepw} informs us that
there exists a connected vertex-transitive cubic graph $X$ of
girth $g$. If $X$ is bipartite, we can set $L_g=X$. Otherwise,
for $L_g$ we can take its bipartite double (the direct product
$X\times K_2$ with the complete graph $K_2$),
which is connected, cubic, bipartite, vertex-transitive, and
has girth $g$. Note that $L_g$ is cyclically $g$-edge-connected
(by \cite[Theorem~17]{NS-cc}) and is also $3$-edge-colourable
because it is bipartite. In order to supply a superedge for our
superposition we take four copies of the $5$-pole $M_g$, append
twenty new vertices to them, and connect all in the manner
illustrated in the left-hand side of
Figure~\ref{fig:superpos}.

To finish the construction, pick a $6$-cycle $C=(e_0e_1\ldots
e_5)$ in $Pg$, and let $v_i$ denote the common vertex of the
edges $e_{i-1}$ and $e_i$ with indices taken modulo $6$.
Substitute the vertices $v_0$ and $v_3$ with a copy of the
$(3,3,1)$-pole $Z$ and the vertices $v_1$, $v_2$, $v_4$, and
$v_5$ with a copy of $M_g$. Next, substitute the edges $e_0$,
$e_2$, $e_3$, and $e_5$ with a copy of the $(3,3)$-pole~$F_g$.
Finally, join the connectors of each copy of $F_g$ to a connector of
a copy of $M_g$ and a connector of a copy of $Z$, and connect
the copies of $M_g$ between themselves and to the remaining
vertices of $Pg$ in such a way that a cubic graph
arises, see Figure~\ref{fig:superpos}.

So far we have not specified the distribution of semiedges of
$M_g$ into connectors that turn a $5$-pole into a
$(3,3,1)$-pole. Nevertheless, it is useful to realise that
whatever the distribution is, $G_g$ is a cyclically
$5$-edge-connected snark of girth $g$. Indeed, the choice of
$L_g$ and the way how superposition has been performed
guarantee that the required cyclic connectivity and girth of
$G_g$ are achieved. 

The proof that $G_g$ is not
$3$-edge-colourable uses the fact that every $3$-edge-colouring
of $F_g$ (that is, a nowhere-zero
$\mathbb{Z}_2\times\mathbb{Z}_2$-flow on $F_g$) induces a
non-zero total flow through $F_g$. The graph $G_g$ thus arises
from $Pg$ by a so-called proper superposition and
\cite[Theorem~4]{Ko}, which implies that $G_g$ is a snark. What
remains is to arrange that the specified pair $\{u,v\}$ of
adjacent vertices of $G_g$  is non-removable. Recall that we
have required the graph $L_g$ to be bipartite and therefore
$3$-edge-colourable. Hence, we can choose a $3$-edge-colouring
of $M_g$ and subsequently, depending on the chosen colouring,
attach the five semiedges  of $M_g$ to the rest of the
superedge $F_g$ in such a way that $F_g$ has a
$3$-edge-colouring that extends to a $3$-edge-colouring of
$G_g-\{u,v\}$. Thus $G_g$ has all the required properties.

\medskip

Now we are ready to prove Theorem~\ref{thm:main}.

\begin{proof}[Proof of Theorem~\ref{thm:main}]
Consider the graph $G_g$ from Theorem~\ref{thm:largegirth} and
let $u$ and $v$ be adjacent vertices such that $G_g-\{u,v\}$ is
$3$-edge-colourable.  Fix an even integer $g\ge 8$.  We show
that for $d=g/2$ the graph $\bar G_g$ obtained from $G_g$ by
inflating both $u$ and $v$ into triangles fulfils the statement
of the theorem. First, we claim that $\df{\bar G_g}=4$. Indeed,
Theorem~\ref{thm:2triangles} implies that $3\le \df{\bar
G_g}\le 4$. However, $\bar G_g$ contains no $6$-cycle so
$\df{\bar G_g}\ne 3$ by Proposition~\ref{prop:core}, whence
$\df{\bar G_g}=4$. To prove that $\rdf{\bar G_g}\ge d$, pick a
regularly optimal $3$-array for $\bar G_g$ and let $K$ be any
component of its core. By the definition, $K$ is a circuit of
an even length, say $k$, that alternates uncovered and doubly
covered edges. It is obvious that $k\ge g$ because every
circuit in $\bar G_g$ is either a triangle or has length at
least $g$. Proposition~\ref{prop:rdf-girth} now yields that
$\rdf{\bar G_g}\ge  k/2\ge g/2=d$,
and the proof is complete.
\end{proof}

\section{Remarks}
\noindent{}The proof of Theorem~\ref{thm:main} demonstrates
that inflating two adjacent vertices of a snark to triangles
may decrease defect from an arbitrary large value to no more
than $4$. By Theorem~\ref{thm:2triangles}, it is sufficient to
inflate a  pair whose removal produces a $3$-edge-colourable
graph. A similar effect can be achieved by inflating a single
vertex if it lies on a \emph{non-removable $5$-cycle}, one
whose removal yields a $3$-edge-colourable graph. Indeed, in
\cite[Theorem~5.3]{KMNS-red} we have constructed a family of
snarks $H_g$, where $g$ is an arbitrary even integer not
smaller than $6$, such that $\df{H_g}\ge g/2$ but after
inflating a vertex on a non-removable $5$-cycle defect
decreases to $3$. Proposition~\ref{prop:rdf-girth} implies that
the same holds for regular defect.

Inflating a vertex of a snark to a triangle may also lead to an
opposite effect. Take a nontrivial snark $G$ of girth at least
$6$ with defect $3$; Isaacs flower snarks $J_n$, where $n$ is
odd, are typical candidates for $G$
(see~\cite[Example~5.4]{KMNS-red}). Choose a subset $D$ of
vertices such that $G-D$ is acyclic and inflate every vertex in
$D$ into a triangle giving rise to a graph $G'$. Every cycle of
length greater than $3$ in $G'$, which corresponds to a cycle
in $G$, has length greater than $6$. Assuming that $G'$ admits
a regular $3$-array we finally get that $\rdf{G'}\geq
\df{G'}\geq 4$.

\medskip

The construction presented in Theorem~\ref{thm:largegirth}
produces snarks with triangles. Since snarks containing
circuits of length smaller than $5$ and cuts of size smaller
than $4$ are often considered to be trivial, the following
problem suggests itself.

\begin{problem}
Does there exist a family of cyclically $4$-edge-connected
snarks with girth at least $5$ such that their defect and
regular defect are arbitrarily far apart?
\end{problem}

\section*{Acknowledgements}
\noindent{}Research reported in this paper was supported by the grant
APVV-23-0076 of Slovak Research and Development Agency. The
first and the third author were partially supported by the
grant VEGA~2/0056/25 of Slovak Ministry of Education. The
second and the fourth authors were partially supported by the grant VEGA 1/0727/22.

The authors are grateful to anonymous referees for their constructive remarks.

\bigskip

\end{document}